\documentclass[oneside, 12pt]{amsart}
\usepackage{amsmath}
\usepackage{amssymb}
\usepackage{color}
\usepackage[usenames,dvipsnames,x11names,svgnames]{xcolor}
\usepackage[colorlinks=true,linkcolor=NavyBlue,citecolor=DarkGreen]{hyperref}
\usepackage{amsthm}
\usepackage{amscd}
\usepackage{geometry}
\usepackage{mathrsfs}

\newtheorem{thm}{Theorem}
\newtheorem*{ethm}{Theorem}
\newtheorem{conj}{Conjecture}
\newtheorem{prop}{Proposition}

\newtheorem{cor}[prop]{Corollary}

\theoremstyle{definition}

\newcommand{\mb}{\mathbb}
\newcommand{\mc}{\mathcal}

\newcommand{\ol}{\overline}
\newcommand{\leqs}{\leqslant }
\newcommand{\geqs}{\geqslant }

\begin{document}
\title[Moments of random multiplicative functions]{Moments of random multiplicative functions and truncated characteristic polynomials}
\author[W. Heap]{Winston Heap}
\address{Department of Mathematical Sciences, Norwegian University of Science and Technology,
NO-7491 Trondheim, Norway} \email{winstonheap@gmail.com}

\author[S. Lindqvist]{Sofia Lindqvist}
\address{Department of Mathematical Sciences, Norwegian University of Science and Technology,
NO-7491 Trondheim, Norway} \email{lindqvist.sofia@gmail.com}

\thanks{Research supported by grant 227768 of the Research Council of Norway.}

\maketitle

\begin{abstract}We give an asymptotic formula for the $2k$th moment of a sum of multiplicative Steinhaus variables. This was recently computed independently by Harper, Nikeghbali and Radziwi\l\l. We also compute the $2k$th moment of a truncated characteristic polynomial of a unitary matrix. This provides an asymptotic equivalence with the moments of Steinhaus variables. Similar results for multiplicative Rademacher variables are given. 
  
\end{abstract}

\section{Introduction}

In the study of the Riemann zeta function there are two probabilistic heuristics which have had significant recent attention. One of these is the use of random multiplicative functions in problems of an arithmetic nature and the other is the use of random matrix theory to model various statistics of the zeta function.  

The study of random multiplicative functions was initiated by Wintner \cite{W} when he modelled the M\"obius function as the multiplicative extension to the squarefree integers of the random variables $\epsilon_p$, each of which takes the values $\{\pm 1\}$ with equal probability. This provided a model for the reciprocal of the Riemann zeta function and hence an appropriate\footnote{We say `appropriate' here since previous models simply used a random $\pm 1$ as the coefficients which, as objected to by Levy \cite{Levy}, did not take into account the multiplicative nature of the problem.} probabilistic interpretation of the Riemann hypothesis. More recently, random models have been used in association with Dirichlet characters in the work of Granville and Soundararajan (e.g. in \cite{GS1, GS2}) and also for the quantities $p^{it}$ when $p$ ranges over the set of primes  \cite{H0, L, LLR}.

The connection between the Riemann zeta function and random matrix theory is well known and has been extensively studied. 
One of the more remarkable predictions of random matrix theory is the Keating--Snaith conjecture \cite{KS} regarding the moments of the zeta function. This states that
\begin{equation}\frac{1}{T}\int_0^T|\zeta(\tfrac{1}{2}+it)|^{2k}dt\sim a(k)g(k)(\log T)^{k^2}\end{equation}
where $g(k)$ is a certain geometric factor involving the Barnes $G$-function and 
 \begin{equation}\label{arithmetic factor}a(k)=\prod_p \bigg(1-\frac{1}{p}\bigg)^{k^2}\sum_{m=0}^\infty \frac{d_k(p^m)^2}{p^m}\end{equation}
with $d_k(n)$ being the $k$-fold divisor function. In essence, the reasoning behind the Keating--Snaith conjecture can be stated as follows.  Since the zeros of the zeta function are conjectured to share the same distribution as eigenvalues of a random matrix  in the CUE, it is reasonable to expect that the characteristic polynomial of a matrix provides a good model to the zeta function in the mean. Thus,  for an appropriate choice of $N$ one could expect
\begin{equation}\label{KS conj 2}\frac{1}{T}\int_0^T|\zeta(\tfrac{1}{2}+it)|^{2k}dt\sim a(k)\mb{E}_{U(N)}\Big[|\Lambda(M,z)|^{2k}\Big]\end{equation}
where $\Lambda(M,z)$ denotes the characteristic polynomial of the matrix $M$ and the expectation is over all unitary matrices of size $N$ with respect to the Haar measure.

Recently, Conrey and Gamburd \cite{CG} showed that the asymptotic in \eqref{KS conj 2} holds if one both truncates the characteristic polynomial and replaces the zeta function by a Dirichlet polynomial of length $x=o(T^{1/k})$. This allowed them to deduce that 
\[\lim_{T\to\infty}\frac{1}{T}\int_0^T \Big|\sum_{n\leq x}n^{-1/2-it}\Big|^{2k}dt\sim a(k)c(k)(\log x)^{k^2},  \qquad k\in\mb{N}\]
where $a(k)$ is given by \eqref{arithmetic factor} and $c(k)$ is the volume of a particular polytope in $\mb{R}^{k^2}$. 
This result may be stated in the more general framework of random multiplicative functions as follows.

Given the set of primes, associate a set of i.i.d. random variables $\{X_p\}$, equidistributed on the unit circle with variance 1. We extend these to the positive integers by requiring that $X_n$ is multiplicative; that is, if $n=\prod_p p^{\alpha(p)}$ then $X_n=\prod_p X_p^{\alpha(p)}$. We let $\mb{E}[\cdot]$ denote the expectation. We refer to the $X_n$ as multiplicative Steinhaus variables. The association $p^{it} \leftrightarrow X_p$ is then seen to be more than just formal in light of the identity
\begin{equation*}\qquad\mb{E}\bigg[\Big|\sum_{n\leq x}X_n/n^\sigma\Big|^{2q}\bigg]=\lim_{T\to\infty}\frac{1}{T}\int_0^T \Big|\sum_{n\leq x}n^{-\sigma-it}\Big|^{2q}dt
\end{equation*}
which holds for all $\sigma\in\mathbb{R}$ and $q>0$. This can be proved by first demonstrating it for $q\in\mathbb{N}$ and then applying the Weierstrass approximation theorem to the function $f:y\mapsto y^{q/2}$. 

Our main aim is to extend the results of Conrey--Gamburd to more general $\sigma$, in particular to $\sigma=0$, and to exhibit the connection between moments of random multiplicative functions and random matrix theory.

\begin{thm}\label{2kth norm thm}For fixed $k\in\mb{N}$ and $0\leq\sigma< 1/2$ we have 

\begin{equation}\label{2kth norm}\qquad\mb{E}\bigg[\Big|\sum_{n\leq x}X_n/n^\sigma\Big|^{2k}\bigg]\sim \frac{a(k)\beta(k)}{(1-2\sigma)^{2k-1}}\frac{\Gamma(2k-1)}{\Gamma(k)^{2}}x^{k(1-2\sigma)}(\log x)^{(k-1)^2}.
\end{equation}
where, $a(k)$ is given by \eqref{arithmetic factor}, $\beta(1)=1$, and 
\begin{equation}\label{beta}
\begin{split}\beta(k)=&\frac{1}{(2\pi i)^{2k-1}}\int_{(b_{2k})}\cdots\int_{(b_2)}
 \bigg[\prod_{i=2}^k\prod_{j=k+1}^{2k}\frac{1}{s_i+s_j}\bigg]e^{s_2+\cdots +s_{2k}}\prod_{j=2}^kds_j\prod_{j={k+1}}^{2k} \frac{ds_j}{s_j}
\end{split}
\end{equation}
for $k\geq 2$. Here, $\int_{(b_j)}=\int_{b_j-i\infty}^{b_j+i\infty}$ and $b_j>0$ for all $j$. 
\end{thm}

In the significant case $\sigma=0$, Theorem \ref{2kth norm thm} has been proved independently by Harper, Nikeghbali and Radziwi\l\l{} \cite{HNR}. It is of interest to note that the constant in their result involves the volume of the Birkhoff polytope $\mathcal{B}_k$. By comparing coefficients we get that $\mathrm{vol}(\mc{B}_k)=k^{k-1}\beta(k)$. It is an open problem to determine a closed form for the volume of the Birkhoff polytope \cite{P} and  a representation in terms of such contour integrals may be new. A direct proof of the equation $\mathrm{vol}(\mc{B}_k)=k^{k-1}\beta(k)$ can be found by applying the methods of section \ref{asymptotics un} to the formula for the Ehrhart polynomial given in \cite{BP}.  Also, our methods work equally well in the case $\sigma=1/2$ and thus by comparing coefficients with Theorem 2 of Conrey--Gamburd \cite{CG}, we get a contour integral representation for their constant (see equation \eqref{alpha}).  

It should be noted that the expectation on the line $\sigma=0$ counts the number of solutions $(m_j)\in\mb{N}^{2k}$ to the equation  $m_1m_2\cdots m_k=m_{k+1}\cdots m_{2k}$ with the restriction $1\leqs m_j\leqs x$. In the case $k=2$, Ayyad, Cochrane and Zheng \cite{ACZ} computed this quantity to a high accuracy. Theorem \ref{2kth norm thm} therefore extends these results to $k\geqs 3$,  although we do not achieve their level of accuracy. By including the extra condition $(m_j,q)=1$ in the  equation, a slight modification of our methods give the following asymptotic formula for Dirichlet character sums.
\begin{thm}\label{char thm}Let $\chi$ be a primitive Dirichlet character modulo $q$ and suppose $q$ has a bounded number of prime factors. Then for fixed $k\in\mb{N}$,
\begin{equation}\label{char formula}\frac{1}{\varphi(q)}\sum_{\chi\neq \chi_0} \Big|\sum_{n\leq x}\chi(n)\Big|^{2k}\sim a(k)\beta(k)\prod_{p|q} \bigg(\sum_{m=0}^\infty\frac{d_k(p^m)^2}{p^m}\bigg)^{-1}\frac{\Gamma(2k-1)}{\Gamma(k)^{2}}x^k(\log x)^{(k-1)^2}
\end{equation} 
as $x, \,q\to\infty$ with $x^k\leq q$ where $\chi_0$ denotes the principal character.
\end{thm}


Our results on the random matrix theory side are as follows. Let $U(N)$ denote the group of unitary matrices of size $N$. For a matrix $M\in U(N)$ with eigenvalues $(e^{i\theta_j})_{j=1}^N$ let 
\[\Lambda(z)=\Lambda(M,z)=\det(I-zM)=\prod_{j=1}^N(1-e^{i\theta_j}z)=\sum_{n=0}^N c_M(n)(-z)^n\]
where the $c_M(n)$ are the secular coefficients.
For $N^\prime\leqs N$ we consider the truncated characteristic polynomial given by
\[\Lambda_{N^\prime}(z)=\sum_{n\leq N^\prime}c_M(n)(-z)^n.\]
Let $\mb{E}_{U(N)}[\,\cdot\,]$ denote the expectation over $U(N)$ with respect to Haar measure.

\begin{thm}\label{rmt thm}Let $k\in\mathbb{N}$ be fixed and suppose $|z|> 1$. Then for all $N\geqs k\log x$ we have
\begin{equation}\mathbb{E}_{U(N)}\big[|\Lambda_{\log x}(z)|^{2k}\big]\sim \frac{\beta(k)}{(1-|z|^{-2})^{2k-1}}\frac{\Gamma(2k-1)}{\Gamma(k)^{2}}F_k(z)x^{2k\log|z|}(\log x)^{(k-1)^2}
\end{equation}
where $\beta(k)$ is given by \eqref{beta} and 
\begin{equation}
\begin{split}\label{F_k}F_k(z)=_2\!F_1(1-k,1-k;2-2k;1-|z|^{-2})
\end{split}\end{equation}
with $_2F_1$ being Gauss' hypergeometric function.
\end{thm}


One may notice a certain similarity between Theorems \ref{2kth norm thm} and \ref{rmt thm}. Indeed, by including the work of \cite{CG} in the case $\sigma=1/2$ we have the following.
\begin{cor}\label{comparison cor}Let $k\in\mb{N}$ be fixed and let $z_\sigma$ be any comlpex number such that $|z_\sigma|=e^{1/2-\sigma}$. Then for $0\leq \sigma\leq1/2$ and $N\geqs k\log x$ we have  
\begin{equation}\label{comparison}\mb{E}\bigg[\Big|\sum_{n\leq x}X_n/n^\sigma\Big|^{2k}\bigg]\sim a(k)c_\sigma(k)\mathbb{E}_{U(N)}\big[|\Lambda_{\log x}(z_\sigma)|^{2k}\big]
\end{equation}
where $a(k)$ is given by \eqref{arithmetic factor} and
\[c_\sigma(k)=\begin{cases}\bigg(\frac{1-e^{2\sigma-1}}{1-2\sigma}\bigg)^{2k-1}F_k(e^{1/2-\sigma})^{-1}, &0\leq \sigma<1/2\\\qquad1, &\sigma=1/2.\end{cases}\]
\end{cor}

A problem which has garnered some attention recently  is to determine the first  moment of $\big|\sum_{n\leqs x} X_n\big|$. A conjecture of Helson \cite{Helson} states that this is $o(\sqrt{x})$, but this seems doubtful now given the evidence in \cite{BS,HNR}. Another motivation for the present article was to provide a conjecture for the first moment via Corollary \ref{comparison cor}. 

Let us then assume that Corollary \ref{comparison cor} holds for $0\leqs k<1$. Then (presumably) the average on the right side of \eqref{comparison} can be taken over matrices of size $N=\log x$ which leads to a computation of the full characteristic polynomial. By an application of Szeg\H o's Theorem, Chris Hughes has shown (\cite{HT}, formula (3.177)) that for $|z|<1$
\begin{equation}\label{Hughes formula}\mathbb{E}_{U(N)}\big[|\Lambda(z)|^{2s}\big]\sim \bigg(\frac{1}{1-|z|^2}\bigg)^{s^2}\end{equation}
as $N\to\infty$. On applying the functional equation 
\[\Lambda(M,z)=\det M (-z)^N \Lambda (M^\dagger,1/z)\]
 with $|z|=e^{\sigma-1/2}$ we obtain the following conjecture.

\begin{conj}\label{steinhaus conj}For $0\leqs k<1$ and $0\leqs \sigma<1/2$ we have 
\begin{equation}\mb{E}\bigg[\Big|\sum_{n\leq x}X_n/n^\sigma\Big|^{2k}\bigg]\sim \frac{a(k)F_k(e^{1/2-\sigma})^{-1}}{(1-e^{2\sigma-1})^{(k-1)^2}(1-2\sigma)^{2k-1}}\,x^{k(1-2\sigma)}.
\end{equation}
\end{conj}

For $k=1/2$ and $\sigma=0$ we can compute the constants to a reasonable accuracy. The arithmetic factor $a(k)$ admits a continuation to real values of $k$ via the formula 
\[d_k(p^m)=\binom{k+m-1}{m}=\frac{\Gamma(k+m)}{m!\Gamma(k)}.\]
We then find that $a(1/2)=0.98849...\,$. The other constants are given by
\[F_{1/2}(e^{1/2})^{-1}=_2\!\!F_1(\tfrac{1}{2},\tfrac{1}{2};1;1-e^{-1})^{-1}=\mathrm{agm}\Big(1-\sqrt{1-\tfrac{1}{e}},1+\sqrt{1-\tfrac{1}{e}}\Big)=0.79099...\]
where agm$(x,y)$ is Gauss' arithmetic-geometric mean and
\[\bigg(\frac{e}{e-1}\bigg)^{1/4}= 1.21250... \,\,.\]
Thus, on combining the constants we acquire the conjecture
\begin{equation}\label{helson conj}\mb{E}\bigg[\Big|\sum_{n\leq x}X_n\Big|\bigg]\sim 0.8769...\sqrt{x}.\end{equation}

One can instead consider multiplicative Rademacher variables. In this case, associate a set of i.i.d. random variables $\{Y_p\}$, which are $\pm 1 $ with uniform probability, to the set of primes. Extend this to all positive integers by requiring $Y_n$ to be multiplicative and non-zero only on the square free integers; that is, $Y_n = \left|\mu(n)\right|\prod_{p|n} Y_p$.

Let 
\begin{equation}\label{b}
  b(k) = \prod_p \left( 1-\frac{1}{p} \right)^{k(2k-1)} \sum_{i=0}^k {2k\choose 2i}\frac{1}{p^i}  
\end{equation}
and
\begin{equation}\label{gamma constant}
  \gamma(k) = \frac{1}{(2\pi i)^{2k}} \int_{(b_{2k})}\cdots\int_{(b_{1})}\prod_{1\le i<j\le 2k}  \frac{1}{s_i+s_j}\prod_{j=1}^{2k} e^{2s_j}ds_j,
\end{equation}
where $b_j>0$ for all $j$. We then have the following result.
\begin{thm}\label{2kth norm thm rad}For fixed $k\in\mb{N}, k\geqs 2$ we have 
\begin{equation}
  \mathbb{E}\left[\biggl|\sum_{n\le x} Y_n\biggr|^{2k}\right]\sim \gamma(k)b(k)2^{2k} x^{k}(\log x)^{2k^2-3k}.
\end{equation}
\end{thm}

Let $SO(2N)$ denote the group of orthogonal $2N\times 2N$ matrices with determinant 1, and let $\mathbb{E}_{SO(2N)}[\cdot]$ denote the expectation over $SO(2N)$ with respect to Haar measure.
\begin{thm}\label{rmt thm so}Let $k\in\mathbb{N}$ be fixed and suppose $z\in \mathbb{R}$, $|z|> 1$. Then for all $N\geqs k\log x$ we have
  \begin{equation}\mathbb{E}_{SO(2N)}\big[|\Lambda_{\log x}(z)|^{2k}\big]\sim \frac{\gamma(k)}{(1-|z|^{-1})^{2k}} x^{2k\log|z|}(\log x)^{2k^2-3k}
\end{equation}
where $\gamma(k)$ is given by \eqref{gamma constant}.
\end{thm}

\begin{cor}\label{comparison rad cor} For fixed $k\in\mb{N}$, $k\geqs 2$ and all $N\geqs k\log x$ we have
  \[\mb{E}\bigg[\Big|\sum_{n\leq x}Y_n\Big|^{2k}\bigg]\sim b(k)2^{2k} (1-e^{-1/2})^{2k}\mathbb{E}_{SO(2N)}\big[|\Lambda_{\log x}(e^{1/2})|^{2k}\big]\]
  where $b(k)$ is the arithmetic factor given by \eqref{b}.
\end{cor}

 Similarly to the case of Steinhaus variables, we expect that the 1st moment is $\sim c\sqrt{x}$ for some constant $c$. Unfortunately we have not been able to find an analogue of \eqref{Hughes formula} for the special orthogonal group and so cannot make a precise conjecture. For some recent results on the order of $\sum_{n\leq x}Y_n$ see \cite{H,LTW}.

\section{Asymptotics for Steinhaus variables: Proof of Theorem \ref{2kth norm thm} }
\subsection{A contour integral representation for the expectation}
We have 
\[\mb{E}\bigg[\Big|\sum_{n\leq x}X_n/n^\sigma\Big|^{2k}\bigg]=\sum_{\substack{n_1\cdots n_k=\\n_{k+1}\cdots n_{2k}\\n_j\leq x}}\frac{1}{(n_1\cdots n_{2k})^\sigma}.\]
We  invoke the condition $n_j\leq x$ in each $j$ by using the contour integral 
\begin{equation}\label{fundamental integral}\qquad\qquad\frac{1}{2\pi i}\int_{(b)}y^{s}\frac{ds}{s}=\begin{cases}1,\,\,&y>1\\ 0,\,\,&y<1\end{cases},\qquad\qquad (b>0)\end{equation}
with $y=x/n_j$. For each $j$ we take a specific line of integration $b_j$. For reasons that will become clear we take $b_1=\epsilon<1-2\sigma$ if $\sigma<1/2$ and $b_1=2$ if $\sigma=1/2$. In both cases we may take the other lines to be sufficiently large so as to guarantee absolute convergence; $b_j=2$ say ($j=2,\ldots, 2k$). This gives
\begin{equation*}
\begin{split}\!\!\!\!\!\!\mb{E}\bigg[\Big|\sum_{n\leq x}X_n/n^\sigma\Big|^{2k}\bigg]=&\sum_{\substack{n_1\cdots n_k=\\n_{k+1}\cdots n_{2k}}}\frac{1}{(n_1\cdots n_{2k})^\sigma}\frac{1}{(2\pi i)^{2k}}\int_{(b_{2k})}\cdots\int_{(b_1)}\prod_{j=1}^{2k} \bigg(\frac{x}{n_j}\bigg)^{s_j}\frac{ds_j}{s_j}\\
=&\frac{1}{(2\pi i)^{2k}}\int_{(b_{2k})}\cdots\int_{(b_1)}F_k(\sigma+s_1,\cdots,\sigma+s_{2k})\prod_{j=1}^{2k} x^{s_j}\frac{ds_j}{s_j}.
\end{split}
\end{equation*}
where
\[F_k(z_1,\cdots,z_{2k})=\sum_{\substack{n_1\cdots n_k=\\n_{k+1}\cdots n_{2k}}}\frac{1}{n_1^{z_1}\cdots n_{2k}^{z_{2k}}}.\]
Since the condition $n_1\cdots n_k=n_{k+1}\cdots n_{2k}$ is multiplicative we may express $F_k(z)$ as an Euler product:
\begin{equation}\label{euler prod}
\begin{split}F_k(z_1,\ldots,z_{2k})=&\prod_p \sum_{\substack{m_1+\cdots +m_k\\=m_{k+1}+\cdots+ m_{2k}}}\frac{1}{p^{m_1z_1+\cdots +m_{2k}z_{2k}}}\\
=&\prod_p \bigg(1+\sum_{i=1}^k\sum_{j=k+1}^{2k}\frac{1}{p^{z_i+z_j}}+O\Big(\sum \frac{1}{p^{z_{i_1}+z_{j_1}+z_{i_2}+z_{j_2}}}\Big)\bigg)\\
=&A_k(z_1,\ldots,z_{2k})\prod_{i=1}^k\prod_{j=k+1}^{2k}\zeta(z_i+z_j) 
\end{split}
\end{equation}
where 
\begin{equation}\label{A_k}A_k(z_1,\ldots,z_{2k})=\prod_p \bigg[\prod_{i=1}^k\prod_{j=k+1}^{2k}\bigg(1-\frac{1}{p^{z_i+z_j}}\bigg)\bigg]\cdot \sum_{\substack{m_1+\cdots +m_k\\=m_{k+1}+\cdots+ m_{2k}}}\frac{1}{p^{m_1z_1+\cdots +m_{2k}z_{2k}}}.\end{equation}
Upon expanding the inner products and sum whilst referring to the middle line of \eqref{euler prod}, we see that $A_k(z_1,\ldots,z_{2k})$ is an absolutely convergent product provided $\Re(z_i+z_j)>1/2$ for $1\leq i\leq k, k+1\leq j\leq 2k$. 

We now have
\begin{multline}\label{norm contour}\mb{E}\bigg[\Big|\sum_{n\leq x}X_n/n^\sigma\Big|^{2k}\bigg]=
\frac{1}{(2\pi i)^{2k}}\int_{(b_{2k})}\cdots\int_{(b_1)}A_k(\sigma+s_1,\ldots,\sigma+s_{2k})\\
\times\prod_{i=1}^k\prod_{j=k+1}^{2k}\zeta(2\sigma+s_i+s_j) \prod_{j=1}^{2k} x^{s_j}\frac{ds_j}{s_j}.
\end{multline}

\subsection{The case $\sigma=1/2$}Although the case $\sigma=1/2$ has already been investigated by Conrey-Gamburd \cite{CG}, we will go over it with a proof that is instructive for the case  $0\leq\sigma<1/2$.

Set $\sigma=1/2$ in \eqref{norm contour}. We write the resulting integral as 
\[\frac{1}{(2\pi i)^{2k}}\int_{(b_{2k})}\ldots\int_{(b_1)}B_k(s_1,\ldots,s_{2k})\prod_{i=1}^k\prod_{j=k+1}^{2k}\frac{1}{s_i+s_j} \prod_{j=1}^{2k} e^{\mc{L}s_j}\frac{ds_j}{s_j}\]
where $\mc{L}=\log x$ and
\[B_k(s_1,\ldots,s_{2k})=A_k(\tfrac{1}{2}+s_1,\ldots,\tfrac{1}{2}+s_{2k})\prod_{i=1}^k\prod_{j=k+1}^{2k}(s_i+s_j)\zeta(1+s_i+s_j).\]
This function is holomorphic in a neighbourhood of $(0,0,\ldots,0)$ and the constant term in its Taylor expansion about this point is given by $A_k(\tfrac{1}{2},\ldots,\tfrac{1}{2})$. 

We now make the substitution $s_j\mapsto s_j/\mc{L}$ in each variable to give an integral of the form
 \[\frac{\mc{L}^{k^2}}{(2\pi i)^{2k}}\int_{(c_{2k})}\cdots\int_{(c_1)}B_k(s_1/\mc{L},\ldots,s_{2k}/\mc{L})\prod_{i=1}^k\prod_{j=k+1}^{2k}\frac{1}{s_i+s_j} \prod_{j=1}^{2k} e^{s_j}\frac{ds_j}{s_j}.\]
First, note that we may shift the contours so as to be independent of $\mc{L}$, to $\Re(s_j)=2$ say. We now truncate the integrals at height $T=o(\mc{L})$ and take a Taylor approximation to $B_k(\underline{s})$ about the point $(0,0,\ldots,0)$. Then upon letting $\mc{L}\to\infty$ we see that this integral is asymptotic to
 \[A_k(\tfrac{1}{2},\ldots,\tfrac{1}{2})\frac{\mc{L}^{k^2}}{(2\pi i)^{2k}}\int_{(b_{2k})}\cdots\int_{(b_1)}\prod_{i=1}^k\prod_{j=k+1}^{2k}\frac{1}{s_i+s_j} \prod_{j=1}^{2k} e^{s_j}\frac{ds_j}{s_j}.\]
A short calculation gives $A_k(\tfrac{1}{2},\ldots,\tfrac{1}{2})=a(k)$ where $a(k)$ is given by \eqref{arithmetic factor}. The remaining constant is given by
\begin{equation}\label{alpha}\alpha(k):=\frac{1}{(2\pi i)^{2k}}\int_{(b_{2k})}\cdots\int_{(b_1)}\prod_{i=1}^k\prod_{j=k+1}^{2k}\frac{1}{s_i+s_j} \prod_{j=1}^{2k} e^{s_j}\frac{ds_j}{s_j}.\end{equation}

We may express $\alpha(k)$ as a volume integral and hence recover the constant of Theorem 2 in \cite{CG}. This is achieved by first writing $(s_i+s_j)^{-1}=\int_0^\infty e^{-x_{ij}(s_i+s_j)}dx_{ij}$ for each term in the product over $i,j$ so that the full product is then given by a $k^2$-fold integral. Upon exchanging the orders of integration and applying \eqref{fundamental integral} the result follows.

\subsection{The case $0\leq\sigma<1/2$}\label{sigma<1/2 section}
Returning to our expression for the expectation given in \eqref{norm contour}, we first make the substitutions $s_j\mapsto s_j+1-2\sigma$ for $k+1\leq j\leq 2k$. This gives
\begin{multline}\label{steinhaus first integral}\mb{E}\bigg[\Big|\sum_{n\leq x}X_n/n^\sigma\Big|^{2k}\bigg]=x^{k(1-2\sigma)}\bigg(\frac{1}{2\pi i}\bigg)^{2k}\int_{(c_{2k})}\cdots\int_{(b_1)}\\
\times A_k(\sigma+s_1,\ldots,\sigma+s_{k},1-\sigma+s_{k+1},\cdots,1-\sigma+s_{2k})\\
\times\prod_{i=1}^k\prod_{j=k+1}^{2k}\zeta(1+s_i+s_j) \prod_{j=1}^{k} x^{s_j}\frac{ds_j}{s_j}\prod_{j={k+1}}^{2k} x^{s_j}\frac{ds_j}{s_j+1-2\sigma}
\end{multline}
In the case $\sigma=1/2$, the leading order term was essentially given by the poles at $s_j=0$. In the present case we must first make the appropriate substitutions to bring the leading order contributions to $s_j=0$. Only then can we make the substitution $s_j\mapsto s_j/\mc{L}$.  

We first extract the polar behaviour of the integrand. Write
\begin{multline}\label{G}G_{k,\sigma}(s_1,\ldots,s_{2k})=A_k(\sigma+s_1,\ldots,\sigma+s_{k},1-\sigma+s_{k+1},\ldots,1-\sigma+s_{2k})\\
\times\prod_{i=1}^k\prod_{j=k+1}^{2k}(s_i+s_j) \zeta(1+s_i+s_j)
\end{multline}
so that our integral becomes
\begin{multline*}
x^{k(1-2\sigma)}\frac{1}{(2\pi i)^{2k}}\int_{(c_{2k})}\cdots\int_{(b_1)}G_{k,\sigma}(s_1,\ldots,s_{2k})
\prod_{i=1}^k\prod_{j=k+1}^{2k}\frac{1}{s_i+s_j}\times\\
\times e^{\mc{L}(s_1+\cdots +s_{2k})}\prod_{j=1}^{k}\frac{ds_j}{s_j}\prod_{j={k+1}}^{2k}\frac{ds_j}{s_j+1-2\sigma}.
\end{multline*}
The function $G_{k,\sigma}(s_1,\ldots,s_{2k})$ is analytic in the region $\Re(s_i+s_j)>-1/2$ for $1\leq i\leq k, k+1\leq j\leq 2k$. 

We now make the substitutions $s_j\mapsto s_j-s_1$ for $k+1\leq j\leq 2k$ and $s_i\mapsto s_i+s_1$ for $2\leq i\leq k$. This gives an integral of the form

\begin{multline*}
x^{k(1-2\sigma)}\frac{1}{(2\pi i)^{2k}}\int_{(d_{2k})}\cdots\int_{(b_1)}G_{k,\sigma}(s_1,s_2+s_1\ldots,s_{k}+s_1,s_{k+1}-s_1,\ldots,s_{2k}-s_1)
\times\\
\times\bigg[\prod_{i=2}^k\prod_{j=k+1}^{2k}\frac{1}{s_i+s_j}\bigg]e^{\mc{L}(s_2+\cdots +s_{2k})} \frac{ds_1}{s_1}\prod_{j=2}^k\frac{ds_j}{s_j+s_1}\prod_{j={k+1}}^{2k} \frac{ds_j}{s_j(s_j-s_1+1-2\sigma)}.
\end{multline*}
Now, for $j=2,3,\ldots,2k$ we let $s_j\mapsto s_j/\mc{L}$. This gives the integral
\begin{multline*}
x^{k(1-2\sigma)}\mc{L}^{(k-1)^2}\frac{1}{(2\pi i)^{2k}}\int_{(e_{2k})}\cdots\int_{(b_1)}\\G_{k,\sigma}(s_1,s_1+s_{2}/\mc{L},\ldots,s_1+s_{k}/\mc{L},-s_1+s_{k+1}/\mc{L},\ldots,-s_1+s_{2k}/\mc{L})
\times\\
\times\bigg[\prod_{i=2}^k\prod_{j=k+1}^{2k}\frac{1}{s_i+s_j}\bigg]e^{s_2+\cdots +s_{2k}} \frac{ds_1}{s_1}\prod_{j=2}^k\frac{ds_j}{\frac{s_j}{\mc{L}}+s_1}\prod_{j={k+1}}^{2k} \frac{ds_j}{s_j(\frac{s_j}{\mc{L}}-s_1+1-2\sigma)}.
\end{multline*}

Once again, we may shift the lines of integration in the integrals over $s_2,s_3,\ldots,s_{2k}$ so as to be independent of $\mc{L}$; back to $\Re(s_j)=2$ say, and truncate the integrals at some height $T=o(\mc{L})$.  From the definition of $G_{k,\sigma}$ given in \eqref{G}, we see that 
\begin{equation*}
\begin{split}\lim_{\mc{L}\to\infty}&G_{k,\sigma}(s_1,s_1+s_{2}/\mc{L},\ldots,s_1+s_{k}/\mc{L},-s_1+s_{k+1}/\mc{L},\ldots,-s_1+s_{2k}/\mc{L})\\
=&A_k(\sigma+s_1,\ldots,\sigma+s_{1},1-\sigma-s_{1},\ldots,1-\sigma-s_{1})\\
=&A_k(0,\ldots,0,1,\ldots,1)\\
=&A_k(\tfrac{1}{2},\ldots,\tfrac{1}{2})
\end{split}
\end{equation*}
where in the last two lines we have used the symmetry of $A_k$. As previously claimed, this last quantity is given by \eqref{arithmetic factor}. The other limits are easily evaluated. 

Thus, as $\mc{L}\to\infty$ we have
\begin{multline}\mb{E}\bigg[\Big|\sum_{n\leq x}X_n/n^\sigma\Big|^{2k}\bigg]\sim a(k)x^{k(1-2\sigma)}\mc{L}^{(k-1)^2}\frac{1}{(2\pi i)^{2k}}\int_{(b_{2k})}\cdots\int_{(b_1)}\\
 \bigg[\prod_{i=2}^k\prod_{j=k+1}^{2k}\frac{1}{s_i+s_j}\bigg]e^{s_2+\cdots +s_{2k}} \frac{ds_1}{s_1^k(1-2\sigma-s_1)^k}\prod_{j=2}^kds_j\prod_{j={k+1}}^{2k} \frac{ds_j}{s_j}.
\end{multline}

For the integral over $s_1$ we push the line of integration to the far left encountering a pole at $s_1=0$. The integral over the new line vanishes in the limit and so 
\begin{equation*}
\begin{split}
\frac{1}{2\pi i}\int_{(b_1)}\frac{ds_1}{s_1^k(1-2\sigma-s_1)^k}=&\frac{1}{(k-1)!}\frac{d^{k-1}}{ds_1^{k-1}}\bigg((1-2\sigma-s_1)^{-k}\bigg)\bigg|_{s_1=0}\\
=&\frac{1}{(k-1)!}\frac{\Gamma\big(k+(k-1)\big)/\Gamma(k)}{(1-2\sigma)^{2k-1}}\\
=&\frac{\Gamma(2k-1)}{\Gamma(k)^2}\frac{1}{(1-2\sigma)^{2k-1}}.
\end{split}
\end{equation*}
The remaining integrals are given by $\beta(k)$ of equation \eqref{beta}.

Note that although it appears as if one should be able to make the substitution $s_j\mapsto s_j/\mc{L}$ directly in \eqref{steinhaus first integral} without first shifting the variables by $s_1$, this is not the case. Upon truncating the integrals at height $T=o(\mc{L})$, the largest error terms arise from the $\zeta$-factors in $G_{k,\sigma}$ when they are evaluated close to $t=0$. For this to occur in all terms of the form $\zeta(1+(s_i+s_j)/\mc{L})$ and $\zeta(1+s_j/\mc{L})$ for $i=2,\dots,k$, $j=k+1,\dots,2k$, one must have $t_i\approx -t_j \approx 0$ for $i=2,\dots,2k$ and $j=k+1,\dots,2k$. When looking at the error arising from cutting some $s_i$ at height $T$, this large contribution clearly is excluded, as one has $|t_i|\ge T\gg 0$ for this $i$. On the other hand, if one makes the substitution $s_j\mapsto s_j/\mc{L}$ directly in \eqref{steinhaus first integral} and attempts to cut all integrals at height $T=o(\mc{L})$, a large error arises from $t_i\approx-t_j$ for $i=1,\dots,k$ and $j=k+1,\dots,2k$.

\section{Character sums: sketch proof of Theorem \ref{char thm}}
We shall only sketch the proof of Theorem \ref{char thm} since it is very similar to the proof of Theorem \ref{2kth norm thm}. Recall the orthogonality property of Dirichlet characters: for $m,n$ coprime to $q$
\[\frac{1}{\varphi(q)}\sum_{\chi }\chi(m)\ol{\chi(n)}=\begin{cases}1\qquad \mathrm{if} \,\,m\equiv n\mod q,\\ 0\qquad\mathrm{otherwise.}\end{cases}\]
This implies that for all $x^k\leq q$
\begin{equation*}
\begin{split}
\frac{1}{\varphi(q)}\sum_\chi \Big|\sum_{m\leq x}\chi(m)\Big|^{2k}=&\frac{1}{\varphi(q)}\sum_\chi\sum_{\substack{m_i\leqs x}}\chi(m_1\cdots m_k)\ol{\chi}(m_{k+1}\cdots m_{2k})\\
=&\sum_{\substack{m_1\cdots m_k=m_{k+1}\cdots m_{2k}\\m_i\leqs x\\(m_i,q)=1}}1.
\end{split}
\end{equation*} 
On applying the line integral \eqref{fundamental integral} we acquire
\[\sum_{\substack{m_1\cdots m_k=m_{k+1}\cdots m_{2k}\\m_i\leqs x\\(m_i,q)=1}}1=\frac{1}{(2\pi i)^{2k}}\int_{(b_{2k})}\cdots \int_{(b_1)} H_{k,q}(s_1,\ldots, s_{2k})\prod_{j=1}^{2k}x^{s_j}\frac{ds_j}{s_j}\]
where 
\[H_{k,q}(s_1,\ldots,s_{2k})=\sum_{\substack{m_1\cdots m_k=m_{k+1}\cdots m_{2k}\\(m_j,q)=1}}\frac{1}{m_1^{s_1}\cdots m_{2k}^{s_{2k}}}.\]
Expressing this as an Euler product gives
\begin{equation*}
\begin{split}
H_{k,q}(s_1,\ldots,s_{2k})=&\prod_{p|q}\bigg(\sum_{\substack{m_1+\cdots+ m_k=\\m_{k+1}+\cdots +m_{2k}}}\frac{1}{p^{m_1{s_1}+\cdots+ m_{2k}{s_{2k}}}}\bigg)^{-1}\prod_p \sum_{\substack{m_1+\cdots+ m_k=\\m_{k+1}+\cdots +m_{2k}}}\frac{1}{p^{m_1{s_1}+\cdots+ m_{2k}{s_{2k}}}}\\
=&\prod_{p|q}\bigg(\sum_{\substack{m_1+\cdots+ m_k=\\m_{k+1}+\cdots +m_{2k}}}\frac{1}{p^{m_1{s_1}+\cdots+ m_{2k}{s_{2k}}}}\bigg)^{-1} A_k(s_1\ldots,s_{2k})\prod_{i,j}\zeta(s_i+s_j)\\
=&C_{k,q}(s_1,\ldots,s_{2k})\prod_{i,j}\zeta(s_i+s_j),
\end{split}
\end{equation*}
say. Here, the function $A_k(s_1,\ldots s_{2k})$ is that of equation \eqref{A_k}. Since the number of prime factors of $q$ remains fixed, $C_{k,q}(s_1,\ldots,s_{2k})$ is holomorphic in the same regions as $A_k(s_1,\ldots,s_{2k})$. The arguments of the previous section now follow, with the arithmetic constant being given by
\begin{equation*}
\begin{split}
C_{k,q}(\tfrac{1}{2},\ldots,\tfrac{1}{2})=&A_k(\tfrac{1}{2},\ldots,\tfrac{1}{2})\prod_{p|q}\bigg(\sum_{\substack{m_1+\cdots+ m_k=\\m_{k+1}+\cdots +m_{2k}}}\frac{1}{p^{m_1+\cdots+ m_{k}}}\bigg)^{-1}\\
=&a(k)\prod_{p|q}\bigg(\sum_{n=0}^\infty\frac{d_k(p^m)^2}{p^{n}}\bigg)^{-1}.
\end{split}
\end{equation*}

Now
\[\frac{1}{\varphi(q)}\Big|\sum_{n\leq x}\chi_0(n)\Big|^{2k}=\frac{1}{\varphi(q)}\Big|\sum_{\substack{n\leq x\\(n,q)=1}}1\Big|^{2k}=\frac{1}{\varphi(q)}\Big(\frac{\varphi(q)}{q}x+O(2^{\omega(q)})\Big)^{2k}\]
where $\omega(q)$ represents the number of distinct prime factors of $q$. Since we're assuming $\omega(q)$ is bounded this last error term is $O(1)$ as $q\to\infty$. Hence,
\[\frac{1}{\varphi(q)}\Big|\sum_{n\leq x}\chi_0(n)\Big|^{2k}\sim \bigg(\frac{\varphi(q)}{q}\bigg)^{2k-1}\frac{x^{2k}}{q}\leq \bigg(\frac{\varphi(q)}{q}\bigg)^{2k-1}x^k.\]
Since this is of a lower order than the main term when $\omega(q)$ is bounded the result follows.

\section{Moments of the truncated characteristic polynomial in the unitary case: Proof of Theorem \ref{rmt thm}}

\subsection{A formula for the expectation}

We begin by recalling the definitions.  Let $U(N)$ denote the group of unitary matrices of size $N$. For a matrix $M\in U(N)$ with eigenvalues $(e^{i\theta_j})_{j=1}^N$ let 
\[\Lambda(z)=\det(I-zM)=\prod_{j=1}^N(1-e^{i\theta_j}z)=\sum_{n=0}^N c_M(n)(-z)^n.\]
The coefficients $c_M(n)$ are known as the secular coefficients. We have $c_M(0)=1$, $c_M(1)=\mathrm{Tr}(M)$ and $c_M(N)=\det(M)$. In general, note that these coefficients are symmetric functions of the eigenvalues.
For $N^\prime\leq N$, consider the truncated characteristic polynomial given by
\[\Lambda_{N^\prime}(z)=\sum_{n\leq N^\prime}c_M(n)(-z)^n.\]
We will compute the expectation of this object as the following multiple contour integral.
\begin{prop}\label{exp cont int prop}Let $k\in\mb{N}$. Then for all $z\in\mb{C}$ and $N\geq k\mc{L}$ we have
 \begin{equation*}\label{exp cont int}\mathbb{E}_{U(N)}\big[|\Lambda_{\mc{L}}(z)|^{2k}\big]=\frac{1}{(2\pi i)^{2k}}\int\cdots\int \frac{(u_1\cdots u_{2k})^{-\mc{L}}}{\prod_{i=1}^k\prod_{j=k+1}^{2k}(1-|z|^2u_iu_j)}\prod_{j=1}^{2k}\frac{du_j}{u_j(1-u_j)}
\end{equation*}
where the integration is around small circles of radii less than $\min(|z|^{-1},1)$. 
\end{prop}

Our plan is to  expand $|\Lambda(z)|^{2k}$, push the expectation through, and then use the results of Diaconis-Gamburd \cite{DG} regarding the expectation of products of the coefficients $c_M(j)$. To state their result we must first detail some notation.

For an $m\times n$ matrix $A$ denote the row and column sums by $r_i$ and $c_j$ respectively and define the vectors
\[\mathrm{row}(A)=(r_1,\ldots,r_m),\]
\[\,\mathrm{col}(A)=(c_1,\ldots,c_n).\]
Given two partitions $\mu=(\mu_1,\ldots,\mu_m)$ and $\tilde{\mu}=(\tilde{\mu}_1,\ldots,\tilde{\mu}_n)$ we let $N_{\mu\tilde{\mu}}$ denote the number of matrices $A$ with $\mathrm{row}(A)=\mu$ and $\mathrm{col}(A)=\tilde{\mu}$. The notation $\langle 1^{a_1}2^{a_2}\cdots\rangle$ is used to represent a partition with $a_1$ parts equal to 1, $a_2$ parts equal to 2 etc. For example, $(5,3,3,2,1)=\langle1^12^13^24^05^1\rangle$.

\begin{ethm}[\cite{DG}]Let $(a_j)_{j=1}^l$, $(b_j)_{j=1}^l$ be sequences of nonnegative integers. Then for $N\geqslant\max\big(\sum_{j=1}^lja_j,\sum_{j=1}^ljb_j\big)$ we have
\[\mb{E}_{U(N)}\bigg[\prod_{j=1}^l c_M(j)^{a_j}\overline{c_M}(j)^{b_j}\bigg]=N_{\mu\tilde{\mu}}\]
where $\mu=\langle 1^{a_1}2^{a_2}\cdots\rangle$ and $\tilde{\mu}=\langle 1^{b_1}2^{b_2}\cdots\rangle$.
\end{ethm}

On expanding the polynomial and pushing the expectation through we get
\begin{equation*}
\begin{split}
&\mathbb{E}_{U(N)}\big[|\Lambda_{\mc{L}}(z)|^{2k}\big]\\=&\sum_{n_1,\ldots,n_{2k}\leq \mc{L}} \mathbb{E}_{U(N)}\big[c_M(n_1)\cdots c_M(n_k)\overline{c_M}(n_{k+1})\cdots \overline{c_M}(n_{2k})\big]z^{n_1+\cdots +n_{k}}\ol{z}^{n_{k+1}+\cdots+n_{2k}}\\
=&\mathbb{E}_{U(N)}\big[|c_M(\mc{L})|^{2k}\big]|z|^{2k\mc{L}}\\&+|z|^{(2k-2)\mc{L}}\sum_{\substack{m\leqs \mc{L}\\n\leqs \mc{L}-1}}\mathbb{E}_{U(N)}\big[c_M(\mc{L})^{k-1}c_M(m)\ol{c_M}(\mc{L})^{k-1}\ol{c_M}(n)\big]z^m\ol{z}^n+\cdots.
\end{split}
\end{equation*}
On taking $N\geqslant k\mc{L}$ the condition of the Theorem is satisfied for all terms in the sum, and can thus be applied. 

Consider the first term. Note that we may write
\begin{equation*}
\begin{split}
\mathbb{E}_{U(N)}\big[|c_M(\mc{L})|^{2k}\big]|z|^{2k\mc{L}}=&\sum_{(m_{ij})_{i,j=1}^k\in B_k(\mc{L})}z^{(\sum_{i=1}^k\sum_{j=1}^k m_{ij})}\ol{z}^{(\sum_{j=1}^k\sum_{i=1}^k m_{ij})}\\
=&\sum_{(m_{ij})_{i,j=1}^k\in B_k(\mc{L})}|z|^{2\sum_{i=1}^k\sum_{j=1}^k m_{ij}}
\end{split}
\end{equation*}
where
\[B_k(\mc{L})=\bigg\{(m_{ij})\in\mathbb{Z}_{\geqslant 0}^{k^2}:\sum_{j=1}^k m_{ij}=\mc{L};\,\sum_{i=1}^k m_{ij}=\mc{L}\bigg\}.\]
Similarly, for the second term we may write
\begin{equation*}
  \begin{split}&|z|^{(2k-2)\mc{L}}\sum_{\substack{m\leqs \mc{L}\\n\leqs \mc{L}-1}}\mathbb{E}_{U(N)}\big[c_M(\mc{L})^{k-1}c_M(m)\ol{c_M}(\mc{L})^{k-1}\ol{c_M}(n)\big]z^m\ol{z}^n\\
=&\sum_{\substack{m\leqs \mc{L}\\n\leqs \mc{L}-1}}\sum_{(m_{ij})_{i,j=1}^k\in C_k(\mc{L},m,n)}z^{\sum_{i=1}^k\sum_{j=1}^k m_{ij}}\ol{z}^{\sum_{i=1}^k\sum_{j=1}^k m_{ij}}
\end{split}
\end{equation*}
where 
\begin{multline*}C_k(\mc{L},m,n)=\bigg\{(m_{ij})\in\mathbb{Z}_{\geqslant 0}^{k^2}:\sum_{j=1}^k m_{ij}=\mc{L}\, \text{ for }\, 0\leqs i\leqs k-1, \sum_{j=1}^km_{kj}=m;\\\,\sum_{i=1}^k m_{ij}=\mc{L}\,\text{  for  }\, 0\leqslant j\leqslant k-1, \sum_{i=1}^k m_{ik}=n\bigg\}.\end{multline*}
This set is empty unless $m=n$ in which case we can write the second term as a sum of $|z|^2$.
Continuing in this fashion we see that
\begin{equation}\mathbb{E}_{U(N)}\big[|\Lambda_{\mc{L}}(z)|^{2k}\big]=\sum_{(m_{ij})_{i,j=1}^k\in D_k(\mc{L})}|z|^{2\sum_{i=1}^k\sum_{j=1}^k m_{ij}}
\end{equation}
where
\[D_k(\mc{L})=\bigg\{(m_{ij})\in\mathbb{Z}_{\geqslant 0}^{k^2}:\sum_{j=1}^k m_{ij}\leqs\mc{L};\,\sum_{i=1}^k m_{ij}\leqs\mc{L}\bigg\}.\]

We now invoke the conditions $\sum m_{ij}\leqs \mc{L}$ with the formula
\begin{equation}\label{fundamental integral 3}\frac{1}{2\pi i}\int_{|u|=\varepsilon}u^{m-\mc{L}}\frac{du}{u(1-u)}=\begin{cases}1,&\,\,\,m\leqs \mc{L}\\ 0 ,&\,\,\,m> \mc{L}\end{cases}\end{equation}
which follows on expanding $(1-u)^{-1}$ as a geometric series.
 This gives
 \begin{multline*}\mathbb{E}_{U(N)}\big[|\Lambda_{\mc{L}}(z)|^{2k}\big]=\sum_{m_{ij}\geqs 0}|z|^{2\sum_{i=1}^k\sum_{j=1}^k m_{ij}}\frac{1}{(2\pi i )^{2k}}\int_{|u_{2k}|=\varepsilon_{2k}}\!\!\!\!\!\!\!\!\cdots\int_{|u_{1}|=\varepsilon_{1}}\times\\
\times u_1^{m_{11}+m_{12}+\cdots+m_{1k}-\mc{L}}u_2^{m_{21}+m_{22}+\cdots+m_{2k}-\mc{L}}\cdots u_k^{m_{k1}+m_{k2}+\cdots+m_{kk}-\mc{L}}\\
\times u_{k+1}^{m_{11}+m_{21}+\cdots+m_{k1}-\mc{L}}u_{k+2}^{m_{12}+m_{22}+\cdots+m_{k2}-\mc{L}}\cdots u_{2k}^{m_{1k}+m_{2k}+\cdots+m_{kk}-\mc{L}}\prod_{j=1}^{2k}\frac{du_j}{u_j(1-u_j)}.
\end{multline*}
On collecting like powers and computing the geometric series we acquire Proposition \ref{exp cont int prop}.

\subsection{Asymptotics for the multiple contour integral}\label{asymptotics un}
Denote the integral  in Proposition \ref{exp cont int prop} by $I$. Thus, 
\[I=\frac{1}{(2\pi i)^{2k}}\int\cdots\int \frac{(u_1\cdots u_{2k})^{-\mc{L}}}{\prod_{i=1}^k\prod_{j=k+1}^{2k}(1-|z|^2u_iu_j)}\prod_{j=1}^{2k}\frac{du_j}{u_j(1-u_j)}\]
where the contours of integration are   positively oriented circles of radii $|u_j|=\varepsilon_j$. We choose $1/|z|^2<\varepsilon_1<1/|z|$, $\varepsilon_j=\varepsilon_1^{-1}|z|^{-2}\delta_j^{1/\mc{L}}$ $(j=k+1,\ldots,2k)$ and $\varepsilon_j=\varepsilon_1\delta_j^{1/\mc{L}}$ ($j=2,\ldots, k$) with $\delta_j<1/|z|^\mc{L}$. With these choices of $\varepsilon_j$ the conditions of Proposition \ref{exp cont int prop} are satisfied, provided $|z|>1$.  
We will now perform similar manipulations to those in section \ref{sigma<1/2 section}

 First, let $u_j \mapsto u_1^{-1}|z|^{-2}u_j$ for $j=k+1,\ldots,2k$. Then
\begin{multline*}I=\frac{x^{2k\log|z|}}{(2\pi i)^{2k}}\int\cdots\int \frac{(u_1^{1-k}u_2\cdots u_{2k})^{-\mc{L}}}{\prod_{i=1}^k\prod_{j=k+1}^{2k}(1-u_1^{-1}u_iu_j)}\\
\prod_{j=1}^{k}\frac{du_j}{u_j(1-u_j)}\prod_{j=k+1}^{2k}\frac{du_j}{u_j(1-u_1^{-1}|z|^{-2}u_j)}.\end{multline*}
Now let $u_j \mapsto u_1 u_j$ for $j=2,\ldots,k$. After a bit of rearranging we have
\begin{multline*}
I=\frac{x^{2k\log|z|}}{(2\pi i)^{2k}}\int\cdots\int \frac{(u_2\cdots u_{2k})^{-\mc{L}}}{\prod_{i=2}^k\prod_{j=k+1}^{2k}(1-u_iu_j)}\frac{du_1}{u_1(1-u_1)}\\
\prod_{j=2}^{k}\frac{du_j}{u_j(1-u_1u_j)}\prod_{j=k+1}^{2k}\frac{du_j}{u_j(1-u_j)(1-u_1^{-1}|z|^{-2}u_j)}.
\end{multline*}

Consider the integral with respect to $u_2,\ldots,u_{2k}$:
\begin{multline*}
J:=\frac{1}{(2\pi i)^{2k-1}}\int\cdots\int \frac{(u_2\cdots u_{2k})^{-\mc{L}}}{\prod_{i=2}^k\prod_{j=k+1}^{2k}(1-u_iu_j)}\prod_{j=2}^{k}\frac{du_j}{u_j(1-u_1u_j)}\\
\prod_{j=k+1}^{2k}\frac{du_j}{u_j(1-u_j)(1-u_1^{-1}|z|^{-2}u_j)}.
\end{multline*}
The contours of integration are given by circles of radii $\varepsilon_j=\delta_j^{1/\mc{L}}$, $j=2,\ldots,2k$. We may now choose the $\delta_j$ to be small but independent of $\mc{L}$ since this does not alter the value of the integral. Then, writing $J$ in its parametrised form gives
\begin{multline*}
J=\frac{1}{(2\pi)^{2k-1}}\int_{-\pi}^{\pi}\cdots\int_{-\pi}^{\pi} \frac{(\delta_2\cdots \delta_{2k})^{-1}e^{-i\mc{L}(\theta_2+\cdots+\theta_{2k})}}{\prod_{i=2}^k\prod_{j=k+1}^{2k}(1-(\delta_i\delta_j)^{1/\mc{L}}e^{i\theta_i+i\theta_j})}\prod_{j=2}^{k}\frac{d\theta_j}{(1-u_1\delta_j^{1/\mc{L}}e^{i\theta_j})}\\
\prod_{j=k+1}^{2k}\frac{d\theta_j}{(1-\delta_j^{1/\mc{L}}e^{i\theta_j})(1-u_1^{-1}|z|^{-2}\delta_j^{1/\mc{L}}e^{i\theta_j})}.
\end{multline*}
Performing the substitutions $\theta_j\mapsto\theta_j/\mc{L}$ and writing $c_j=\log \delta_j$ we have
\begin{multline*}
J=\frac{\mc{L}^{(k-1)^2}}{(2\pi)^{2k-1}}\int_{-\pi\mc{L}}^{\pi\mc{L}}\cdots\int_{-\pi\mc{L}}^{\pi\mc{L}} \frac{e^{-(c_2+i\theta_2+\cdots+c_{2k}+i\theta_{2k})}}{\prod_{i=2}^k\prod_{j=k+1}^{2k}\mc{L}(1-e^{(c_i+i\theta_i+c_j+i\theta_j)/\mc{L}})}\\
\prod_{j=2}^{k}\frac{d\theta_j}{(1-u_1e^{(c_j+i\theta_j)/\mc{L}})}\prod_{j=k+1}^{2k}\frac{d\theta_j}{\mc{L}(1-e^{(c_j+i\theta_j)/\mc{L}})(1-u_1^{-1}|z|^{-2}e^{(c_j+i\theta_j)/\mc{L}})}.
\end{multline*}
We now divide by $\mc{L}^{(k-1)^2}$ and take the limit as $\mc{L}\to\infty$. By an argument involving dominated convergence we may pass the limit through the integral. Then, since $\mc{L}(1-e^{-z/\mc{L}})\sim z$, we acquire
\begin{multline*}
J\sim\frac{\mc{L}^{(k-1)^2}}{(1-u_1)^{k-1}(1-u_1^{-1}|z|^{-2})^k}\frac{1}{(2\pi)^{2k-1}}\int_{-\infty}^{\infty}\cdots\int_{-\infty}^{\infty} \\\frac{e^{-(c_2+i\theta_2+\cdots+c_{2k}+i\theta_{2k})}}{\prod_{i=2}^k\prod_{j=k+1}^{2k}(- c_i-i\theta_i-c_j-i\theta_j)}
\prod_{j=2}^{k}d\theta_j\prod_{j=k+1}^{2k}\frac{d\theta_j}{(-c_j-i\theta_j)}.
\end{multline*}
Finally, we let $\theta_j\mapsto -\theta_j$. Upon noting that $-c_j$ is postive, we see that this last integral is the parametrised form of $\beta(k)$ of equation \eqref{beta}. Therefore,
\[
I\sim\beta(k)x^{2k\log|z|}\mc{L}^{(k-1)^2}\frac{1}{2\pi i}\int\frac{1}{(1-u_1)^{k}(1-u_1^{-1}|z|^{-2})^k}\frac{du_1}{u_1}.
\]

It remains to show that 
\begin{equation}
\begin{split}
I_1:=&\frac{1}{2\pi i}\int\frac{1}{(1-u_1)^{k}(1-u_1^{-1}|z|^{-2})^k}\frac{du_1}{u_1}=\frac{1}{2\pi i}\int\frac{u_1^{k-1}}{(1-u_1)^k(u_1-|z|^{-2})^k}du_1\\
=&\frac{1}{(1-|z|^{-2})^{2k-1}}\frac{\Gamma(2k-1)}{\Gamma(k)^{2}}F_k(z)
\end{split}
\end{equation}
where \begin{equation}
\begin{split}\label{F_k}F_k(z)=_2\!F_1(1-k,1-k;2-2k;1-|z|^{-2})
\end{split}\end{equation}
and $_2F_1$ is the usual hypergeometric function.

Note that in $I_1$ we are still integrating over a circle of radius $\varepsilon<1/|z|$ since we made no substitutions in the variable $u_1$. Thus, the only contribution is from the pole at $u_1=|z|^{-2}$. Therefore,
\begin{equation}\label{I1}
\begin{split}I_1=&\frac{1}{(k-1)!}\frac{d^{k-1}}{du^{k-1}}\bigg(\frac{u^{k-1}}{(1-u)^k}\bigg)\bigg|_{u=|z|^{-2}}\\
=&\frac{1}{(k-1)!}\sum_{m=0}^{k-1}\binom{k-1}{m}\frac{d^m}{du^m}\big[(1-u)^{-k}\big]\frac{d^{k-1-m}}{du^{k-1-m}}\big[u^{k-1}\big]\bigg|_{u=|z|^{-2}}\\
=&\frac{1}{\Gamma(k)}\sum_{m=0}^{k-1}\binom{k-1}{m}\cdot\frac{\Gamma(k+m)/\Gamma(k)}{(1-|z|^{-2})^{k+m}}\cdot\frac{\Gamma(k)}{\Gamma(m+1)}|z|^{-2m}\\
=&\frac{1}{(1-|z|^{-2})^k}\cdot\frac{1}{\Gamma(k)}\sum_{m=0}^{k-1}(-1)^m\binom{k-1}{m}\frac{\Gamma(k+m)}{\Gamma(m+1)}\bigg(\frac{1}{1-|z|^2}\bigg)^m.
\end{split}
\end{equation}
As a quick aside, we note that
\begin{equation}\label{2F1}\frac{1}{\Gamma(k)}\sum_{m=0}^{k-1}(-1)^m\binom{k-1}{m}\frac{\Gamma(k+m)}{\Gamma(m+1)}\bigg(\frac{1}{1-|z|^2}\bigg)^m={_2F_1}\big(1-k,k;1;1/(1-|z|^2)\big)\end{equation}
and that this can be written in terms of the Legendre polynomials $P_n(x)$ via the formula (see \cite{BE}, section 3.2)
\[_2F_1(-\lambda,\lambda+1;1;z)=P_\lambda(1-2z).\]

 Continuing with our manipulations, the last line of equation \eqref{I1} may be rewritten as  
\begin{multline}\frac{|z|^{2k}}{(|z|^{2}-1)^{2k-1}}\frac{1}{\Gamma(k)}\sum_{m=0}^{k-1}(-1)^m\binom{k-1}{m}\frac{\Gamma(2k-1-m)}{\Gamma(k-m)}({1-|z|^2})^m\\
=\frac{|z|^{2k}}{(|z|^{2}-1)^{2k-1}}\frac{\Gamma(2k-1)}{\Gamma(k)^2}\sum_{m=0}^{k-1}\frac{(1-k)_m}{m!}\frac{(1-k)_m}{(2-2k)_m}({1-|z|^2})^m\end{multline}
where $(x)_m$ is rising factorial or Pochammer symbol defined by
\[(x)_m=\begin{cases}1, \hspace{4.5cm} m=0,\\x(x+1)\cdots(x+m-1), \hspace{0.4cm}m\geqs 1.\end{cases}\]
 This last sum is the hypergeometric function $_2F_1(1-k,1-k;2-2k;1-|z|^2)$. By formula (18) in section 2.9 of \cite{BE} we have
\[ {}_2F_1(1-a,1-b;2-c;z)=(1-z)^{b-1}{}_2F_1(a+1-c,1-b;2-c;z/(z-1)).\] 
On applying this the result follows.

\section{Moments of Rademacher variables: Proof of Theorem \ref{2kth norm thm rad} }
\subsection{A contour integral representation for the norm}
Note that $\mathbb{E}\left[ Y_{n_1}\cdots Y_{n_{2k}} \right] = 1$ if $n_1\cdots n_{2k}$ is a square number and $n_i$ is square-free for $i=1,\dots,2k$, and otherwise it equals zero. Therefore we have
\begin{equation*}
   \mathbb{E}\left[\biggl|\sum_{n\le x} Y_n\biggr|^{2k}\right] = \sum_{\substack{n_1\cdots n_{2k} \text{ square}\\ n_j \le x}}|\mu(n_1)|\cdots |\mu(n_{2k})|.
\end{equation*}
As earlier we invoke the condition $n_j\le x$ in each $j$ by using \eqref{fundamental integral} with $y=x/n_j$.  For each $j$ we integrate along the lines $b_j=2$. This gives
\begin{equation}\label{rademacher integral}
  \begin{split}
    \mathbb{E}\left[\biggl|\sum_{n\le x} Y_n\biggr|^{2k}\right] =& \sum_{\substack{n_1\cdots n_{2k} \text{ square}\\|\mu(n_j)|=1}} \frac{1}{(2\pi i)^{2k}}\int_{(b_{2k})}\cdots \int_{(b_1)} \prod_{j=1}^{2k} \left( \frac{x}{n_j} \right)^{s_j} \frac{ds_j}{s_j}\\
  =& \frac{1}{(2\pi i)^{2k}}\int_{(b_{2k})} \cdots \int_{(b_1)} F_{k}(s_1,\dots,s_{2k}) \prod_{j=1}^{2k} x^{s_j} \frac{ds_j}{s_j}
\end{split}
\end{equation}
where
\begin{equation*}
  F_{k}(z_1,\dots,z_{2k}) = \sum_{\substack{n_1\cdots n_{2k} \text{ square}}} \frac{|\mu(n_1)|\cdots |\mu(n_{2k})|}{n_1^{s_1}\cdots n_{2k}^{s_{2k}}}.
\end{equation*}
Since the condition $n_1\cdots n_{2k}$ being square is multiplicative we may express $F_{k}(z)$ as an Euler product:
\begin{equation*}
  \begin{split}
    F_{k} (z_1,\dots,z_{2k}) =& \prod_p \sum_{\substack{m_1+\cdots + m_{2k} \text{ even}\\0\leqs m_j\leqs 1}} \frac{1}{p^{m_1z_1+\cdots m_{2k} z_{2k}}}\\
    =& \prod_p \left( 1+\sum_{1\le i<j\le 2k} \frac{1}{p^{z_i+z_j}} + O\left( \sum \frac{1}{p^{z_{i_1}+\cdots+z_{i_4}}} \right) \right)\\
    =& B_{k}(z_1,\dots,z_{2k}) \prod_{1\le i<j\le 2k} \zeta(z_i+z_j)
  \end{split}
\end{equation*}
where
\begin{equation*}
  B_{k}(z_1,\dots,z_{2k}) = \prod_p \left[ \prod_{1\le i<j\le 2k} \left( 1-\frac{1}{p^{z_i+z_j}} \right) \right]\cdot \sum_{\substack{m_1+\cdots + m_{2k} \text{ even}\\0\leqs m_j\leqs 1}} \frac{1}{p^{m_1z_1+\cdots +m_{2k}z_{2k}}}.
\end{equation*}
Similarly as in the proof of Theorem \ref{2kth norm thm} we get that $B_{k}(z)$ is an absolutely convergent product provided $\Re(z_i+z_j)>1/2$ for $1\le i<j\le 2k$.

Now we make the substitution $s_j \mapsto s_j+1/2$ for $j=1,\dots,2k$ in the second line of \eqref{rademacher integral} to get
\begin{equation*}
  \mathbb{E}\left[\biggl|\sum_{n\le x} Y_n\biggr|^{2k}\right] = x^k \frac{1}{(2\pi i)^{2k}}\int_{(b_1^\prime)}\cdots \int_{(b_{2k}^\prime)} G_{k}(s_1,\dots,s_{2k}) \prod_{1\le i<j\le 2k} \frac{1}{s_i+s_j} \prod_{j=1}^{2k} x^{s_j} \frac{ds_j}{s_j +\frac{1}{2}}
\end{equation*}
where we define
\begin{equation*}
  G_{k}(z_1,\dots,z_{2k}) = B_k(s_1+1/2,\dots,s_{2k}+1/2) \prod_{1\le i<j\le 2k} \zeta(1+s_i+s_j)(s_i+s_j).
\end{equation*}
The function $G_{k}(z)$ is analytic in the region $\Re (z_i+z_j)>-1/2$ for $1\le i<j\le 2k$, and $G_{k}(0,\dots,0)=B_{k}(1/2,\dots,1/2)$. Finally we make the substitution $s_j\mapsto s_j/\mc{L}$ for $j=1,\dots ,2k$ to get
\begin{equation*}
    x^{k}\mc{L}^{2k^2-3k} \frac{1}{(2\pi i)^{2k}}\int\cdots\int G_{k}(s_1/\mc{L},\dots,S_{2k}/\mc{L})\prod_{1\le i<j\le 2k} \frac{1}{s_i+s_j}\prod_{j=1}^{2k} e^{s_j} \frac{ds_j}{s_j/\mc{L}+\frac{1}{2}}.
\end{equation*}

Shift the lines of integration to be independent of $\mc{L}$, say back to $\Re(s_j)=2$ for $j=1,\dots,2k$ and truncate each line at height $T=o(\mc{L})$. Computing the Taylor expansions and then taking the limit as $\mc{L}\to\infty$ gives
\begin{equation*}
  \mathbb{E}\left[\biggl|\sum_{n\le x} Y_n\biggr|^{2k}\right] \sim b(k)2^{2k}x^k \mc{L}^{2k^2-3k} \gamma(k)
\end{equation*}
where $b(k)=B_k(1/2,\dots,1/2)$ is the arithmetic factor given in \eqref{b} and $\gamma(k)$ is the integral given by \eqref{gamma constant}.

\section{Moments of the truncated characteristic polynomial in the special orthogonal case: Proof of Theorem \ref{rmt thm so}}

\subsection{A formula for the expectation.} 

As for the $U(N)$ case we begin by expressing the expectation as a multiple contour integral.

\begin{prop}\label{exp cont int prop so}Let $k\in\mb{N}$, $\mc{L}>1$. Then for all $z\in\mb{R}$ and $N\geq k\mc{L}$ we have
 \begin{equation*}\label{exp cont int so}\mathbb{E}_{SO(2N)}\big[\Lambda_{\mc{L}}(z)^{2k}\big]=\frac{1}{(2\pi i)^{2k}}\int\cdots\int \frac{(u_1\cdots u_{2k})^{-\mc{L}}}{\prod_{1\le i<j\le 2k}(1-z^2u_iu_j)}\prod_{j=1}^{2k}\frac{du_j}{u_j(1-u_j)}
  \end{equation*}
  where the integration is around small circles of radii less than $\min(|z|^{-1},1)$.
\end{prop}
To prove this we must use an alternative method to before since there is no analogous result of Diaconis--Gamburd \cite{DG} for the special orthogonal group. Instead, we  use the following result of Conrey-Farmer-Keating-Rubinstein-Snaith.

\begin{ethm}[\cite{CFKRS}]Let $d\mu$ denote the Haar measure on $SO(2N)$. Then for $m\geqs 1$ we have
  \begin{multline*}
    \int_{SO(2N)} \Lambda(w_1)\cdots \Lambda(w_{m}) d\mu =\\ w_1^N \cdots w_m^N \left[\sum_{\epsilon_j\in \{1,-1\}} \left(\prod_{j=1}^m w_j^{N\epsilon_j}\right) \prod_{1\le i<j\le m} \frac{1}{1-w_i^{-\epsilon_i}w_j^{-\epsilon_j}}\right].
  \end{multline*}
\end{ethm}
To begin with, we use the integral
\begin{equation}\label{fundamental integral 2}
  \frac{1}{2\pi i}\int_{(c)} e^{Ys} \frac{ds}{s} = \begin{cases} 1,\,\,& Y>0\\ 0,\,\,& Y<0 \end{cases} \qquad\qquad (c>0)
\end{equation}
to write
\begin{equation*}
  \Lambda_\mc{L} (z) = \frac{1}{2\pi i}\int_{(c)} \Lambda(ze^{-s}) e^{\mc{L}s} \frac{ds}{s}
\end{equation*}
for each of the factors $\Lambda_\mc{L}(z)$. Pushing through the expectation then gives
\begin{equation*}
  \mathbb{E}_{SO(2N)}\left[ \Lambda_\mc{L} (z)^{2k} \right] = \frac{1}{(2\pi i)^{2k}}\int_{(c_{2k})} \cdots \int_{(c_1)} \mathbb{E}\left[ \prod_{i = 1}^{2k} \Lambda(ze^{-s_i})\right]\prod_{i=1}^{2k} e^{\mc{L} s_i} \frac{ds_i}{s_i}
\end{equation*}
where, for reasons that will become apparent, we take $c_1 > c_2 > \dots > c_{2k} > \max \{0,\log|z|\}$ and $c_i - c_{i-1} > 2\log|z|$ for $i=2,\dots,2k$. Using the theorem with $w_i=ze^{-s_i}$ for $i=1,\dots,2k$ then gives
\begin{multline*}
  \mathbb{E}_{SO(2N)}\left[ \Lambda_\mc{L} (z)^{2k} \right] = \sum_{\epsilon_j\in\{1,-1\}} \frac{1}{(2\pi i)^{2k}} \int_{(c_{2k})} \cdots \int_{(c_1)} |z|^{2Nk} \left( \prod_{i=1}^{2k} z^{N\epsilon_i} e^{-N\epsilon_i s_i} \right)\times \\ \times \prod_{1\le i<j\le 2k} \frac{1}{1-z^2e^{\epsilon_i s_i + \epsilon_j s_j}} \prod_{j=1}^{2k} e^{(\mc{L}-N)s_j} \frac{ds_j}{s_j}.
\end{multline*} 

We would like to show that any term with $\epsilon_j=1$ for some $j\in\{1,\dots,2k\}$ gives zero contribution. To this end, let $n=\min\{j\in\{1,\dots,2k\}:\epsilon_j=1\}$. When integrating over $s_1,\dots,s_{n-1}$ we only need to keep track of the highest power of $e^{s_n}$. For the $s_1$ integral, write
\begin{equation*}
  \frac{1}{1-z^2e^{-s_1+\epsilon_j s_j}} = \sum_{m\ge 0} z^{2m} e^{m(-s_1+\epsilon_j s_j)}
\end{equation*}
for $j=2,\dots,2k$. Using \eqref{fundamental integral 2} we see that $m\le \mc{L}$ in each sum, so the highest possible contribution of powers of $e^{s_n}$ is $e^{\mc{L}s_n}$. Further, as $\epsilon_2=\cdots=\epsilon_{n-1}$, this integral contributes a nonpositive power of $e^{s_2},\dots,e^{s_{n-1}}$.

Continuing in this fashion and integrating $s_2,\dots,s_{n-1}$, we deduce that one is left with
\begin{multline*} 
  \frac{1}{(2\pi i)^{2k-n+1}} \int_{(c_{2k})} \cdots \int_{(c_{n+1})}\left(\int_{(c_n)} \prod_{j=n+1}^{2k} \frac{1}{1-z^2 e^{s_n +\epsilon_j s_j}}e^{(\mc{L}-2N)s_n}\left( e^{(n-1)\mc{L} s_n} + \cdots \right) \frac{ds_n}{s_n} \right) \times \\
  \times \left(\prod_{j=n+1}^{2k} z^{N\epsilon_j} e^{(\mc{L}-N-N\epsilon_j) s_j}\right)\prod_{n+1\le i<j\le 2k} \frac{1}{1-z^2 e^{\epsilon_i s_i + \epsilon_j s_j}} \prod_{j=n+1}^{2k} \frac{ds_j}{s_j}
\end{multline*}
multiplied by some power of $z$, where the additional terms in $\left( e^{(n-1)\mc{L}s_n}+\cdots \right)$ are all lower powers of $e^{s_n}$. For the innermost integral, expanding the factors in power series gives
\begin{multline*}
  \frac{1}{2\pi i} \int_{(c_n)} \prod_{j=n+1}^{2k} \frac{1}{1-z^2e^{s_n+\epsilon_j s_j}} e^{(\mc{L}-2N)s_n} \left( e^{(n-1)\mc{L}s_n}+\dots \right) \frac{ds_n}{s_n} \\
  = \frac{1}{2\pi i}\int_{(c_n)} \left(\prod_{j=n+1}^{2k} z^{-2} e^{-s_n-\epsilon_j s_j} \sum_{m\ge 0} x^{-2m} e^{(-s_n-\epsilon_j s_j)m}\right) e^{(\mc{L}-2N)s_n}\left( e^{(n-1)\mc{L} s_n}+\dots \right) \frac{ds_n}{s_n}.
\end{multline*}
The highest possible power of $e^{s_n}$ amongst all terms is $e^{(2k-n +(n-1)\mc{L} +\mc{L}-2N)s_n}$ which is negative for $N\ge \mc{L}k$ and $\mc{L} > 1$. By \eqref{fundamental integral 2} all terms are thus zero.

This leaves
\begin{equation*}
  \mathbb{E}_{SO(2N)}\left[ \Lambda_\mc{L} (z)^{2k} \right] = \frac{1}{(2\pi i)^{2k}}\int_{(c_{2k})} \cdots \int_{(c_1)} \prod_{1\le i<j\le 2k} \frac{1}{1-z^2e^{-s_i-s_j}}\prod_{j=1}^{2k} e^{\mc{L}s_j} \frac{ds_j}{s_j}.
\end{equation*}
In order to arrive at the contour integral, we expand each factor of the product as
\begin{equation*}
  \frac{1}{1-z^2 e^{-s_i-s_j}} = \sum_{m_{ij}\ge 0} z^{2m_{ij}} e^{(-s_i-s_j)m_{ij}}.
\end{equation*}
Separating the integrals and using \eqref{fundamental integral 2} then gives
\begin{equation*}
  \mathbb{E}_{SO(2N)}\left[ \Lambda_\mc{L} (z)^{2k} \right] = \sum_{(m_{ij})_{i,j=1}^{2k}\in A_k(\mc{L})} z^{2\sum_{1\le i<j\le 2k} m_{ij}}
\end{equation*}
where
\begin{equation*}
  A_k(\mc{L}) = \bigg\{ (m_{ij})\in \mathbb{Z}_{\ge 0}^{k(2k-1)}:\sum_{j=1}^{i-1} m_{ji}+\sum_{j=i+1}^{2k}m_{ij} \le \mc{L}, i=1,\dots,2k\bigg\}.
\end{equation*}

As for $U(N)$ we invoke the conditions $\sum m_{ij}\le \mc{L}$ with the formula \eqref{fundamental integral 3}, giving
\begin{multline*}
  \mathbb{E}_{SO(2N)}\left[ \Lambda_\mc{L} (z)^{2k} \right] = \sum_{m_{ij}\ge 0} z^{2\sum_{1\le i<j\le 2k}m_{ij}} \frac{1}{(2\pi i)^{2k}}\int_{|u_{2k}|=\varepsilon_{2k}}\cdots\int_{|u_{1}|=\varepsilon_{1}} \times\\
  u_1^{m_{12} + m_{13}+\cdots+m_{1,2k}-\mc{L}}u_2^{m_{12}+m_{23}+\cdots+m_{2,2k}-\mc{L}} \cdots u_{2k}^{m_{1,2k}+m_{2,2k}+\cdots+m_{2k-1,2k}-\mc{L}} \prod_{j=1}^{2k} \frac{du_j}{u_j(1-u_j)}.
\end{multline*}
On collecting like powers and computing the geometric series we aquire Proposition \ref{exp cont int prop so}.

\subsection{Asymptotics for the multiple contour integral}
Denote the integral in Proposition \ref{exp cont int prop so} by $I$. Again we perform manipulations similar to those in section \ref{asymptotics un}.

Let $u_j\mapsto |z|^{-1} u_j$ for $j=1,\dots,2k$. Then
\begin{equation*}
  I = \frac{x^{2k\log |z|}}{(2\pi i)^{2k}}\int\cdots\int \frac{(u_1\cdots u_{2k})^{-\mc{L}}}{\prod_{1\le i<j\le 2k} (1-u_iu_j)} \prod_{j=1}^{2k} \frac{du_j}{u_j(1-|z|^{-1}u_j)}.
\end{equation*}
Next, let $u_j\mapsto u_j^{1/\mc{L}}$ for $j=1,\dots,2k$. This gives
\begin{equation*}
  I = \frac{x^{2k\log|z|}}{(2\pi i)^{2k}}\int\cdots\int \frac{(u_1\cdots u_{2k})^{-1}}{\prod_{1\le i<j\le 2k}(1-(u_i u_j)^{1/\mc{L}})}\prod_{j=1}^{2k} \frac{du_j}{\mc{L}u_j(1-|z|^{-1}u_j^{1/\mc{L}})}
\end{equation*}
which can be expressed as
\begin{equation*}
  \frac{I}{x^{2k\log|z|}\mc{L}^{2k^2-3k}} = \frac{1}{(2\pi i)^{2k}}\int\cdots\int \frac{(u_1\cdots u_{2k})^{-1}}{\prod_{1\le i<j\le 2k} \mc{L}(1-(u_iu_j)^{1/\mc{L}})}\prod_{j=1}^{2k} \frac{du_j}{u_j(1-|z|^{-1}u_j^{1/\mc{L}})}.
\end{equation*}
Here, the contours wind around the origin $\mc{L}$ times. We now choose the radii of the contours to be independent of $\mc{L}$ and then write the integral in parametrised form. Upon taking the limit as $\mc{L}\to \infty$ and pushing the limit through the integrals we acquire
\begin{equation*}
  I \sim \frac{x^{2k\log|z|}\mc{L}^{2k^2-3k}}{(1-|z|^{-1})^{2k}} \gamma(k)
\end{equation*}
where $\gamma(k)$ is given by \eqref{gamma constant}.

\section{Concluding remarks}

Admittedly, our evidence for Conjecture \ref{steinhaus conj} is rather weak and there is a certain level of ambiguity in choosing the size $N$ of the matrices. However, it is interesting that for our choice of $N=\log x$ the random matrix expectation seems to capture the phase change that we expect to see from the expectation of the Steinhaus variables. Indeed, if the conjecture holds for $k=1/2$, then we can obtain the order of magnitude predicted by the conjecture for $0\leqs k\leqs 1$ using H\"olders inequality since we know the value at $k=1$. Also, it seems a little strange, but not impossible, that one could obtain more than square-root cancellation in the case $k=1/2$ as conjectured by Helson.   

Finally, we note the following argument taken from \cite{BS} which gives an upper bound on the Steinhaus expectation for $k=1/2$.
Let $0\leqs u,v<1$ and let $S_x=\sum_{n\leqs x}X_n$. By the Cauchy--Schwarz inequality we have 
\begin{equation}
\begin{split}
\mb{E}\big[|S_x|\big]^2\leqs& \mb{E}\big[|(1-uX_2)(1-vX_3)S_x|^2\big]\cdot \mb{E}\big[ |(1-uX_2)(1-vX_3)|^{-2}\big]\\
=&\frac{1}{(1-u^2)(1-v^2)} \mb{E}\big[|(1-uX_2)(1-vX_3)S_x|^2\big]
\end{split}
\end{equation}
Now,
\begin{equation}
\begin{split}\mb{E}\big[|(1&-uX_2)(1-vX_3)S_x|^2\big]\\
=&\mb{E} \bigg[\Big(1-2u\Re(X_2)+u^2-2v\Re(X_3)+4uv\Re(X_2)\Re(X_3)-2u^2v\Re(X_3)\\
&\,\,+v^2-2uv^2\Re(X_2)+u^2v^2\Big)|S_x|^2\bigg]\\
\sim&x\Big(1-u+u^2-\tfrac{2}{3}v+\tfrac{2}{3}uv-\tfrac{2}{3}u^2v+v^2-uv^2+u^2v^2\Big).
\end{split}
\end{equation}
In this last line we have expanded the square of $S_x$ and used
\[\mb{E}\Big[\sum_{m,n\leq x}X_{am}\ol{X}_{bn}\Big]=\sum_{\substack{am=bn\\m,n\leqs x}}1\sim \frac{1}{ab}x.\]
For $0\leqs u,v<1$, the minimum of the function 
\[f(x,y)=\frac{1-u+u^2-\tfrac{2}{3}v+\tfrac{2}{3}uv-\tfrac{2}{3}u^2v+v^2-uv^2+u^2v^2}{(1-u^2)(1-v^2)}\]
is found to be $\approx0.8164965809...$. Taking square roots gives 
\[\mb{E}\big[|S_x|\big]\leqs (1+o(1))\cdot0.903...\sqrt{x}.\]
Of course, further optimisations may prove to disprove conjecture \eqref{helson conj}.

\section{Acknowledgements}The authors would like to thank Chris Hughes for making us aware of \eqref{Hughes formula}, and also for some useful discussions. We would also like to thank Kristian Seip and Brad Rodgers for their helpful remarks and also  Maksym Radziwi\l\l \, for pointing out the reference \cite{BP}.


\begin{thebibliography}{}



\bibitem{A}G. Andrews, R. Askey, R. Roy, \emph{Special functions}, Cambridge university press, 1999.



\bibitem{ACZ}A. Ayyad, T. Cochrane, Z. Zheng, \emph{The congruence $x_1x_2\equiv x_3x_4$ (mod p), the equation $x_1x_2=x_3x_4$ and the mean value of character sums} J. Number Theory 59 (1996), 398--413.



\bibitem{BE} H. Bateman, A. Erdelyi, \emph{Higher Transcendental Functions Volume I},  McGraw-Hill 1953.


\bibitem{BP}M. Beck, D. Pixton, \emph{The Ehrhart polynomial of the Birkhoff polytope},  Discrete Comput. Geom., {\bf 30} no. 4 (2003) 623--637.


\bibitem{BS}A. Bondarenko, K. Seip, \emph{Helson's problem for sums of a random multiplicative function}, Preprint available at \href{http://arxiv.org/abs/1411.6388}{arXiv.1411.6388}. 


\bibitem{CFKRS} J.B. Conrey, D.W. Farmer, J.P. Keating, M.O. Rubinstein, N.C. Snaith, \emph{Autocorrelation of Random Matrix Polynomials}, Commun. Math. Phys. no. 237 (2003) 365--395.


\bibitem{CG} B. Conrey and A. Gamburd, \emph{Pseudomoments of the Riemann zeta-function and pseudomagic squares}, J. Number Theory \textbf{117} (2006), 263--278.



\bibitem{DG}P. Diaconis, A. Gamburd, \emph{Random matrices, magic squares and matching polynomials}, Electron. J.
Combin. 11 (2) (2004) \#R2.


\bibitem{GS1}A. Granville, K. Soundararajan \emph{Large Character Sums}, J. Am. Math. Soc  {\bf 14} (2001) 365--397.


\bibitem{GS2}A. Granville, K. Soundararajan, \emph{The distribution of values of $L(1,\chi_d)$}, Geom. Funct. Anal. {\bf 13} no. 5 (2003)  992--1028.


\bibitem{H0}A. Harper,  \emph{A note on the maximum of the Riemann zeta function, and log-correlated random variables}. Preprint available at \href{http://arxiv.org/abs/1304.0677}{arxiv.1304.0677}.


\bibitem{H} A. Harper, \emph{Bounds on the suprema of Gaussian processes, and omega results for the sum of a random multiplicative function}, Ann. App. Prob. {\bf 23} no. 2 (2013) 584--616.


\bibitem{HNR}A. Harper, A. Nikeghbali, M. Radziwi\l\l, \emph{A note on Helson's conjecture on moments of random multiplicative functions}, to appear in ``Analytic Number Theory" in honor of Helmut Maier's 60th birthday. Preprint available at \href{http://arxiv.org/abs/1505.01443}{arxiv.1505.01443}.


\bibitem{Helson} H. Helson, \emph{Hankel forms}, Studia Math. {\bf 198} no. 1 (2010)  79--84.


\bibitem{HT} C. P.  Hughes, \emph{On the characteristic polynomial of a random unitary matrix and the Riemann zeta function}, PhD thesis, University of Bristol, 2001.


\bibitem{KS}J. Keating, N. Snaith, \emph{Random matrix theory and $\zeta(1/2+it)$}, Commun. Math. Phys. {\bf 214}  (2000) 57 -- 89.


\bibitem{K} M. S. Klamkin, \emph{Extensions of Dirichlet's multiple integral}, SIAM J. Math. Anal. {\bf 2} no. 3 (1971) 467--469.


\bibitem{L}Y. Lamzouri,  \emph{The two dimensional distribution of values of $\zeta(1+it)$}, Int. Math. Res. Not. IMRN {\bf 2008} Art.106, 48 pp. 


\bibitem{LLR}Y. Lamzouri, S. Lester, M. Radziwi\l\l, \emph{Discrepancy bounds for the distribution of the Riemann zeta-function and applications}, Preprint available at \href{http://arxiv.org/abs/1402.6682}{arxiv.1402.6682}.

    
\bibitem{LTW} Y. Lau, G. Tenenbaum, J. Wu, \emph{Mean values of random multiplicative functions}, Proc. Amer. Math. Soc. {\bf 141} no. 2 (2013) 409--420.


\bibitem{Levy}P. L\'evy, \emph{Sur les s\'eries dont les termes sont des variables eventuelles ind\'ependantes}, Studia
Math. {\b 3 } (1931) 119--155.


\bibitem{P}I. Pak, \emph{Four questions on Birkhoff polytope}, Ann. Comb. {\bf 4} (2000) 83--90.


\bibitem{W}A. Wintner, \emph{Random factorizations and Riemann's hypothesis}, Duke Math. J. {\bf 11} no. 2 (1944) 267--275.



  
\end{thebibliography}
\end{document}